\documentclass[11pt,reqno]{amsart}
\usepackage[T1]{fontenc}
\usepackage[utf8]{inputenc}
\usepackage{lmodern}
\usepackage{amsmath,amssymb,mathtools}
\usepackage{xcolor}
\usepackage{cancel}

\usepackage{microtype}
\usepackage[margin=1in]{geometry}
\usepackage{enumitem}
\usepackage{comment}
\usepackage{graphicx}
\usepackage{needspace}
\usepackage[colorlinks=true,linkcolor=blue,citecolor=blue,urlcolor=blue]{hyperref}
\allowdisplaybreaks[2]
\numberwithin{equation}{section}
\newtheorem{theorem}{Theorem}[section]

\newtheorem{lemma}[theorem]{Lemma}
\theoremstyle{definition}
\newtheorem{definition}[theorem]{Definition}
\theoremstyle{remark}
\newtheorem{remark}[theorem]{Remark}
\newcommand{\Ric}{\operatorname{Ric}}
\newcommand{\Vol}{\operatorname{Vol}}
\newcommand{\tr}{\operatorname{tr}}
\newcommand{\Hess}{\operatorname{Hess}}
\newcommand{\Sph}{\mathbb S}
\newcommand{\R}{\mathbb R}
\newcommand{\dd}{\,d}
\setlist[enumerate,1]{label=\textup{(\roman*)},leftmargin=*,itemsep=3pt}

\usepackage{marginnote}

\title[Distributional scalar curvature and small ball volumes]
{Distributional scalar curvature and small ball volumes}
\author{Changliang Wang}
\address{School of Mathematical Sciences and Institute for Advanced Study,
Key Laboratory of Intelligent Computing and Applications (Ministry of Education),
Tongji University, Shanghai 200092, China}
\email{wangchl@tongji.edu.cn}
\author{Zhixin Wang}
\address{Department of Mathematics, Shanghai Jiao Tong University,
Shanghai 201100, China}
\email{jhin@sjtu.edu.cn}
\date{}

\begin{document}
\begin{abstract}
We prove that distributional scalar curvature lower bounds imply
the small ball volume test for $C^1$ metrics. This is done by
expressing the leading volume correction as a pairing of
distributional scalar curvature with a suitable weight function
and estimating the remainder. We also show that the $C^1$
assumption cannot in general be weakened to continuity, even
under additional Sobolev regularity.
Finally, by constructing examples, we show that
the analogous volume test for fixed symmetric tangent regions
can fail under $C^{1,\alpha}$ convergence, even for smooth
approximating metrics with a uniform positive scalar curvature
lower bound.
\end{abstract}
\maketitle

\tableofcontents

\section{Introduction}\label{sec:introduction}

Various approaches have been proposed to define scalar curvature for
metrics with low regularity and for limit spaces. For continuous
Riemannian metrics with $W^{1,2}_{\mathrm{loc}}$ regularity, the
distributional approach interprets scalar curvature through
integration by parts, building on the work of LeFloch--Mardare
\cite{LeFlochMardare2007} and Lee--LeFloch
\cite{LeeLeFloch2015}. Burkhardt-Guim \cite{BurkhardtGuim2019}
introduced weak pointwise lower scalar curvature bounds for
continuous metrics on closed manifolds using Ricci
flow. Geometric approaches include Gromov's proposals based on
mean convex domains with corners and dihedral angle comparison
\cite{Gromov2014}, as well as his notion of $n$-volumic scalar
curvature, studied further by Deng \cite{Deng2021}, which uses
small ball volume comparisons on metric measure spaces. {Closely
related compactness and stability questions for manifolds with
scalar curvature bounded below and their intrinsic flat limits
are surveyed by Sormani \cite{Sormani-conjectures}.} In the setting of
noncollapsed pointed Gromov--Hausdorff limits with a uniform lower
Ricci curvature bound, Naber \cite[Conjecture~2.20]{Naber2020}
formulated an energy identity conjecture describing the weak
convergence of scalar curvature measures. A basic problem is to
determine how these weak notions are related and under which
regularity and convergence assumptions lower scalar curvature
bounds are preserved.

An important result in this direction is Gromov's $C^0$ stability
theorem \cite{Gromov2014}, for which Bamler gave a Ricci flow proof
\cite{Bamler2016}. If $g_j$ are $C^2$ metrics converging locally
uniformly to a $C^2$ metric $g$, and the scalar curvature $R_{g_j}\ge\kappa$ for a
continuous function $\kappa$, then $R_g\ge\kappa$. The assumption
that the limit is $C^2$ allows its scalar curvature to be interpreted
classically. For less regular limits, one must first specify the
curvature condition to be preserved.

{In this paper, we mainly investigate the relationship between scalar curvature lower bound conditions in the sense of volume limit test (see Definition \ref{def:ball}) and distribution (see Definition \ref{def:distributional-scalar}) for Riemannian metrics of $C^{1}$ regularity.}

For the volume limit test,  we first recall that on a smooth
Riemannian n-dimensional manifold, the volume of a small geodesic ball
has asymptotic expansion
\begin{equation}\label{eq:smooth-ball}
 \Vol_g(B_g(p,t))
 =\omega_n t^n-\frac{\omega_n}{6(n+2)}R_g(p)t^{n+2}
       +O(t^{n+4}),
\end{equation}
where $B_g(p, t)$ is the geodesic ball centered at $p$ with radius $t$, $R_g(p)$ is the scalar curvature at p, and $\omega_n$ is the Euclidean unit ball volume. Small ball volume comparison is one of the approaches
discussed by Gromov in \cite{Gromov2023}. It has the advantage that
it can be formulated on a metric measure space using only distances
and measures of balls, without taking derivatives of the metric. 

\begin{definition}[Volume limit test]\label{def:ball}
 A continuous Riemannian metric $g$ on an $n$-manifold satisfies the
 {\em volume limit test} for scalar curvature bounded below by a constant $K\in\R$ at $p$ if
 \begin{equation}\label{eq:ball-test}
  \limsup_{t\downarrow0}
  \frac{\Vol_g(B_g(p,t))-\omega_n t^n}{t^{n+2}}\le-\frac{\omega_n K}{6(n+2)}.
 \end{equation}
\end{definition}


For a smooth metric, the distance ball volume expansion \eqref{eq:smooth-ball} shows that the condition in Definition \ref{def:ball}
is equivalent to $R_g(p)\ge K$. The behavior of scalar curvature under weak convergence makes it
important to understand which local geometric conditions survive
in the limit. Lee--Topping \cite{LeeTopping2026} showed that, in dimensions $n\ge4$, the distance function of any Riemannian metric conformal to the round metric on $\Sph^n$ is a uniform limit of distance
functions induced by metrics of positive scalar curvature. In dimension three, Kazaras--Xu \cite{KazarasXu2023} constructed three dimensional manifolds with almost nonnegative scalar curvature converging to drawstring spaces in Sormani-Wenger's intrinsic flat topology. These results highlight the distinction between convergence of distances and convergence of metric tensors, and further motivate studying the relation between weak scalar curvature bounds and infinitesimal ball volume inequalities under explicit regularity assumptions.

Another natural approach to defining weak nonnegative scalar curvature uses integration by parts to interpret scalar curvature in the distributional sense \cite{LeeLeFloch2015,LeFlochMardare2007}. We recall the formulation of Lee--LeFloch \cite{LeeLeFloch2015} as the following.

\begin{definition}[Scalar curvature distribution, \cite{LeeLeFloch2015}]\label{def:distributional-scalar}
 Let $g$ be a continuous Riemannian metric on a smooth manifold
$M$, with $g\in W^{1,2}_{\mathrm{loc}}$.
All derivatives below are understood in the weak sense.
 In a coordinate neighborhood $U$, put $J=\sqrt{\det(g_{ij})}$ and
 \[
  \Gamma^k_{ij}
  =\frac12g^{k\ell}
    \bigl(\partial_i g_{j\ell}
          +\partial_j g_{i\ell}
          -\partial_\ell g_{ij}\bigr).
 \]
 Define
 \begin{equation}\label{eq:density-divergence}
  V^k
  =J\left(g^{ab}\Gamma^k_{ab}-g^{ak}\Gamma^b_{ab}\right)
 \end{equation}
 and
 \begin{equation}\label{eq:density-F}
   F={}Jg^{ab}
     \left(\Gamma^k_{k\ell}\Gamma^\ell_{ab}
           -\Gamma^k_{b\ell}\Gamma^\ell_{ak}\right)-\partial_k(Jg^{ab})\Gamma^k_{ab}
       +\partial_b(Jg^{ab})\Gamma^k_{ak},
 \end{equation}
 where repeated indices are summed. Since $g$ and $g^{-1}$
are locally bounded, we have
$V\in L^2_{\mathrm{loc}}$ and $F\in L^1_{\mathrm{loc}}$.
The {\em scalar curvature distribution}, with the volume density
included, is defined locally by
 \begin{equation}\label{eq:local-pairing}
  \langle R_g, u \rangle
  =-\int_U V^k\partial_k u \dd x
    +\int_U F u \dd x,
  \qquad u\in C_c^\infty(U).
 \end{equation}
 For smooth $g$, this convention gives $\langle R_g, u \rangle
 =\int_M R_g u \dd\mu_g.$
 These local expressions define a distribution on $M$ by a
 partition of unity.

 The above definition is a local coordinate form of
Lee--LeFloch's definition, obtained by taking the Euclidean
background metric and expanding the derivative of
$J u$ in their pairing formula. As explained after
Definition~2.1 in \cite[Section~2.1]{LeeLeFloch2015},
the resulting distribution is independent of the background
metric for $g\in C^0\cap W^{1,2}_{\mathrm{loc}}$.
Thus our local expressions agree on overlapping charts
and define the same scalar curvature distribution.

 We say that $g$ has \emph{distributional scalar curvature bounded below by a constant $K\in\R$},
 and write $R_g\ge K$, if
 \begin{equation}\label{eq:distributional-positive}
  \langle R_g, u \rangle\ge K\int_M u \dd\mu_g
  \qquad\text{for every nonnegative } u \in C_c^\infty(M).
 \end{equation}
\end{definition}

Distributional scalar curvature also arises naturally in compactness
and rigidity problems. In various warped product settings, Park--Tian-Wang \cite{ParkTianWang2018},  Tian--Wang \cite{TianWang2024} and the authors \cite{WangWang2026} verified that the limit metrics in Gromov-Sormani MinA scalar compactness conjecture satisfy distributional nonnegative scalar curvature. Positive mass theorems under distributional nonnegative scalar curvature condition were studied and established by  Lee--LeFloch \cite{LeeLeFloch2015}, 
Jiang--Sheng--Zhang \cite{JiangShengZhang2022}, Li \cite{Li2020}, Grant-Tassotti \cite{GrantTassotti2012} and others. For recent developments in regularization and positive mass theorem for
distributional scalar curvature, see e.g. Lee-Litzinger-Simon \cite{LeeLitzingerSimon2026}, Hafemann \cite{Hafemann2026} and Graf-Hafemann \cite{GrafHafemann2026}.



Our main theorem shows that a distributional scalar curvature lower bound $K$ implies the volume limit test scalar curvature lower bound for $C^1$ metrics. It also gives a quantitative upper remainder
under $C^{1}$ regularity.

\begin{theorem}\label{thm:distributional-balls}
 Let $M$ be a smooth $n$-dimensional manifold without boundary,
 $n\ge2$, and let $g$ be a $C^1$ Riemannian metric with
 distributional scalar curvature bounded below by a constant
 $K\in\R$.
 For every $p\in M$, there exist $t_p>0$ and a nonnegative
 function $\varepsilon_p(t)\to0$ as $t\downarrow0$ such that
 \begin{equation}\label{eq:distributional-ball-bound}
  \Vol_g(B_g(p,t))
  \le \omega_n t^n-\frac{\omega_n K}{6(n+2)}t^{n+2}
       +\varepsilon_p(t)t^{n+2},
  \qquad 0<t<t_p.
 \end{equation}
 More precisely, in coordinates centered at $p$ with
 $g(0)=g_{\mathrm E}$ and $Dg(0)=0$, one may take
 \[
  \varepsilon_p(t)
  =C_p\bigl(\rho_p(4t)^2+|K|t\rho_p(4t)\bigr),
  \qquad
  \rho_p(r):=\sup_{|x|\le r}|Dg(x)|.
 \]
 In particular, $g$ satisfies \eqref{eq:ball-test} at every point.

 If $g_i$ are $C^1$ Riemannian metrics converging to $g$
in $C^1_{\mathrm{loc}}$ and each $g_i$ has
distributional scalar curvature bounded below by $K$, then
 $t_p$ and $\varepsilon_p(t)$ in
 \eqref{eq:distributional-ball-bound} can be chosen
 {uniformly} for $g$ and all $g_i$.
\end{theorem}

To prove Theorem \ref{thm:distributional-balls}, we expresses the first variation of the volume of a distance
ball as the pairing of the linearized scalar curvature with a
nonpositive radial weight. This weight may be singular at the center,
but it belongs to $W^{1,1}$ after extension by zero. We justify its
pairing with the scalar curvature distribution by sign preserving
approximation. After normalizing the metric and its first derivatives
at the center, the nonlinear error is quadratic in the $C^1$ size of
the rescaled perturbation. Continuity of the first derivatives makes
this error $o(t^{n+2})$.

The relation between scalar curvature and volume, through the
adjoint of the linearized scalar-curvature operator, also appears
in the deformation results of Corvino--Eichmair--Miao
\cite{CorvinoEichmairMiao2013} and the volume comparisons of Yuan
\cite{Yuan2023}. Here the domain is a distance ball that depends on
the metric, and the error is controlled using only $C^1$ regularity. 

The regularity assumption in
Theorem~\ref{thm:distributional-balls} cannot in general be
replaced by continuity and Sobolev regularity, even when the
distributional scalar curvature is represented by a continuous
nonnegative function.

\begin{theorem}\label{thm:below-c1}
 Set $\beta=\sqrt{3}-1$. There exist a neighborhood $U$ of the
 origin in $\R^3$ and a Riemannian metric $g$ on $U$ such that
 \[
  g\in C^{0,\beta}(U)\cap W^{1,p}_{\mathrm{loc}}(U)
  \qquad\text{for every }1\le p<6+3\sqrt{3}.
 \]
The metric is smooth on $U\setminus\{0\}$, and its
distributional scalar curvature is represented by a
continuous nonnegative function. Moreover, there exists
$c>0$ such that
 \[
  \Vol_g(B_g(0,t))
  =\frac{4\pi}{3}t^3+c\,t^{3+2\beta}
     +O(t^{3+3\beta})
  \qquad\text{as }t\downarrow0.
 \]
 In particular,
 \[
  \lim_{t\downarrow0}
  \frac{\Vol_g(B_g(0,t))-\frac{4\pi}{3}t^3}{t^5}
  =+\infty.
 \]
\end{theorem}

These results address the relation between distributional scalar
curvature and infinitesimal ball volumes. For $C^1$ metrics,
a distributional scalar curvature lower bound implies the
corresponding small ball volume inequality. At lower regularity,
however, this implication can fail even when the scalar curvature
distribution is represented by a continuous nonnegative function.

The volume limit test admits a natural extension from geodesic balls
to more general regions with certain symmetry defined through the exponential map.
Let $\Omega\subset\mathbb{R}^n$ be a bounded measurable set
containing the origin, and assume the following symmetry condition holds:
\begin{equation}\label{omega-condition}
\int_\Omega x^i x^j\,dx=c_\Omega\delta_{ij},
\qquad
c_\Omega:=\frac1n\int_\Omega |x|^2\,dx>0.
\end{equation}
For example, this condition holds if $\Omega$ is invariant under
coordinate reflections and permutations.

Given a smooth Riemannian manifold $(M^n,g)$, identify
$\mathbb{R}^n$ with $T_pM$ using a $g_p$-orthonormal basis and define
\[
\Omega_r(p):=\exp_p(r\Omega).
\]
Integrating the standard volume density expansion in normal
coordinates {(see \cite{Gray1974})} over $r\Omega$, we obtain
\begin{equation}
\begin{aligned}\label{omega-volume-expansion}
\operatorname{Vol}_g(\Omega_r(p))
&=
r^n|\Omega|
-\frac{r^{n+2}}6
\sum_{i,j=1}^n
\operatorname{Ric}_{ij}(p)
\int_\Omega x^i x^j\,dx
+O(r^{n+3})\\
&=
r^n|\Omega|
-\frac{c_\Omega}{6}R_g(p)\,r^{n+2}
+O(r^{n+3}).
\end{aligned}
\end{equation}
Thus the symmetry condition \eqref{omega-condition} cancels the contribution
of the trace-free part of the Ricci tensor, and only
the scalar curvature left. In particular,
\[
\lim_{r\downarrow0}
\frac{r^n|\Omega|-\operatorname{Vol}_g(\Omega_r(p))}
     {r^{n+2}}
=
\frac{c_\Omega}{6}R_g(p).
\]
For convenience, we call the condition that the left hand side is
nonnegative the {\em $\Omega$-volume limit test}. The next result shows that the { $\Omega$-volume limit test }might fail even under $C^{1,\alpha}$ limits.

\begin{theorem}\label{thm:example}
 There is a fixed bounded measurable star-shaped set
 $\Omega\subset\R^3$ satisfying \eqref{omega-condition} with the following
 property. For every $0<\alpha<1$ and every $K>0$, there are smooth
 metrics $g_j$ on $M=\Sph^2\times\Sph^1$ and a $C^{1,\alpha}$ metric
 $g$ such that
 \[
  g_j\longrightarrow g\quad\text{in }C^{1,\alpha},\qquad
  R_{g_j}\ge\frac12.
 \]
 They have the form
 \begin{equation}\label{eq:metric-form}
  g_j=\varphi_j^{-2}h+\varphi_j^2d\xi^2,
  \qquad g=\varphi^{-2}h+\varphi^2d\xi^2,
 \end{equation}
 where $h$ is the unit round metric on $\Sph^2$ and the functions
 $\varphi_j,\varphi$ are positive and independent of $\xi$.
 For some $p\in\Sph^2$, put $\Gamma=\{p\}\times\Sph^1$. Then:
 \begin{enumerate}
  \item $g$ is smooth on $M\setminus\Gamma$ and is the standard
  product metric outside a neighborhood of $\Gamma$;
  \item at each $q\in\Gamma$, under the coordinate orthonormal
  identification $T_qM\simeq\R^3$, for all sufficiently small $t>0$, every initial velocity in
$t\Omega$ determines a unique geodesic on $[0,1]$ starting at $q$,
the resulting exponential map is injective on $t\Omega$
and has a measurable image, and
  \begin{equation}\label{eq:bad-volume}
   \lim_{t\downarrow0}
   \frac{\Vol_g(\exp_q^g(t\Omega))-t^3|\Omega|}{t^5}>K.
  \end{equation}
 \end{enumerate}
\end{theorem}

Metrics of the form \eqref{eq:metric-form} appear in
\cite{KazarasXu2023}. The construction uses the fact that
\begin{equation}\label{eq:scalar-formula-intro}
 R_g=\varphi^2R_h-2|\nabla_h\varphi|_h^2,
\end{equation}
whereas, in coordinates $h=dr^2+\cos^2r\,d\theta^2$,
\begin{equation}\label{eq:ricci-intro}
 \Ric_g(\partial_r,\partial_\theta)
 =-\frac{\varphi_{r\theta}}{\varphi}
  -2\frac{\varphi_r\varphi_\theta}{\varphi^2}
  -\tan r\,\frac{\varphi_\theta}{\varphi}.
\end{equation}
For a smooth metric, the mixed Ricci contribution in
\eqref{omega-volume-expansion} cancels by the symmetry condition
\eqref{omega-condition}. The idea is to break this cancellation.
We choose $\varphi$ so that $\varphi_{r\theta}$ changes sign
between the smooth sectors, and choose $\Omega$ so that the
corresponding geodesics stay within these sectors. The resulting mixed Ricci contributions do not cancel and cause
the $\Omega$-volume limit test to fail.

The paper is organized as follows. Section~\ref{sec:balls}
derives the distance ball volume expansion, expresses its leading
correction term as a scalar curvature pairing, and proves
Theorem~\ref{thm:distributional-balls}.
Section~\ref{sec:below-c1} constructs the counterexample in
Theorem~\ref{thm:below-c1}, showing that the convergence assumption
in Theorem~\ref{thm:distributional-balls} cannot in general be replaced by
$C^0_{\mathrm{loc}}$ convergence.
The final section constructs a counterexample to the
$\Omega$-volume limit test under $C^{1,\alpha}$ convergence, even with
a uniform positive scalar curvature lower bound.

\medskip
{\bf Declaration of AI use.}
ChatGPT was used to assist with the construction of the
counterexample in Theorem~\ref{thm:below-c1} and to improve
the language and presentation of the manuscript.
All other mathematical results and proofs were developed
by the authors. The authors take full responsibility for
the content of the manuscript.

\medskip
{\bf Declaration of conflicts of interest.}
The authors declare that they have no conflicts of interest.

\medskip
{\bf Data availability.}
Data sharing is not applicable to this article, as no datasets
were generated or analysed during this study.
\medskip

\noindent\textbf{Declaration of AI use.}
ChatGPT assistant was used to improve the language
and presentation of the manuscript, and to construct the example in Theorem \ref{thm:below-c1}. All other mathematical
results and proofs are the authors' own.
The authors take full responsibility for the
content of the manuscript.

\medskip
{\bf Acknowledgements.}
 CW is partially supported by the Natural Science Foundation of Shanghai (Grant No. 25ZR1401357) and the Fundamental Research Funds for the Central Universities. 
 ZW is supported by the China Postdoctoral Science Foundation (Grant No. 2025T180838).


\section{Distributional scalar curvature and small distance balls}
\label{sec:balls}

For a continuous Riemannian metric $g$, the notation
$B_g(p,r)$ denotes the open ball of radius $r$ centered at $p$ in
its induced distance. All derivatives, traces, tensor norms, and integration measures
without a metric subscript are with respect to the Euclidean metric $g_{\mathrm E}$. 

The following lemma gives a first order expansion of
$\Vol_g(B_g(0,1))$ for $g=g_{\mathrm E}+h$ with small
$\|h\|_{C^1}$, together with a uniform quadratic remainder estimate.
\begin{lemma}\label{lem:variation}
 Let $g=g_{\mathrm E}+h$ be a $C^1$ Riemannian metric on
 $B_{g_{\mathrm E}}(0,4)\subset\R^n$, where $h$ is a symmetric 2-tensor,
 and suppose $\varepsilon=\|h\|_{C^1(B_{g_{\mathrm E}}(0,4))}$ is sufficiently small. Then
 \begin{equation}\label{eq:volume-variation}
  \Vol_g(B_g(0,1))=\omega_n+\mathcal L_n(h)+O(\varepsilon^2),
 \end{equation}
 where
 \begin{equation}\label{eq:functional}
  \mathcal L_n(h)
  =\frac12\int_{B_{g_{\mathrm E}}(0,1)}\tr h\dd x
   -\frac12\int_{S^{n-1}}\int_0^1
        h(s\theta)(\theta,\theta)\dd s\dd\theta.
 \end{equation}
 The error term in \eqref{eq:volume-variation} is independent of higher
 derivatives of $h$.
\end{lemma}
The first term in $\mathcal L_n(h)$ accounts for the first order
change in the volume density $\sqrt{\det g}$, while the second term
accounts for the difference between $B_{g_{\mathrm E}}(0,1)$ and $B_g(0,1)$.

\begin{proof}
 Throughout the proof, $C$ denotes a universal positive constant that may change from line to line and depends only on $n$ and an upper bound for $\varepsilon$.

 Fix $x \in \mathbb{R}^n$ with $1/2\le|x|\le3/2$. Uniform closeness to the Euclidean
 metric implies that every curve from $0$ reaching
 $\partial B_{g_{\mathrm E}}(0,3)$ is longer than the straight
 segment to $x$. The direct method for the energy on the compact
 ball $\overline{B_{g_{\mathrm E}}(0,3)}$ therefore yields a
 $g$-minimizing curve $\gamma$ from $0$ to $x$ that stays in
 $B_{g_{\mathrm E}}(0,3)$. Parametrize $\gamma$ on $[0,1]$ at
 constant $g$-speed. The Euler--Lagrange equation implies that $\gamma$ is a $C^2$
geodesic, since $g$ is $C^1$. Comparison with the straight segment and
 uniform equivalence of $g$ and $g_{\mathrm E}$ give
 $|\dot\gamma|_g\le C$ and $|\dot\gamma|\le C$.
 Write $\gamma(s)=sx+\eta(s)$, so that $\eta(0)=\eta(1)=0$.
 The geodesic equation $\ddot\gamma^k
 =-\Gamma^k_{ij}(\gamma)\dot\gamma^i\dot\gamma^j$,
 together with the bound $|\Gamma_g|\le C\varepsilon$ following
 from the $C^1$ bound on $h$, gives
 \[
  \|\ddot\eta\|_{L^\infty}
  =\|\ddot\gamma\|_{L^\infty}\le C\varepsilon.
 \]
 Moreover, using $\eta(0)=\eta(1)=0$, we have
 \begin{equation}\label{dis-estimate-eta}
  \|\eta\|_{L^\infty}
  +\|\dot\eta\|_{L^\infty}
  +\|\ddot\eta\|_{L^\infty}
  \le C\varepsilon.
 \end{equation}

 Using $\dot\gamma=x+\dot\eta$, the triangle inequality, and
 \eqref{dis-estimate-eta}, we obtain
 \begin{align}
  &\left|\int_0^1h(\gamma)(\dot\gamma,\dot\gamma)\dd s
        -\int_0^1h(sx)(x,x)\dd s\right| \notag\\
  &\quad\le
  \int_0^1\left|\bigl(h(\gamma(s))-h(sx)\bigr)(x,x)\right|\dd s
  +2\int_0^1\left|h(\gamma(s))(x,\dot\eta(s))\right|\dd s
  \notag\\
  &\qquad\quad
  +\int_0^1
    \left|h(\gamma(s))(\dot\eta(s),\dot\eta(s))\right|\dd s
  \notag\\
  &\quad\le
  |x|^2\|Dh\|_\infty\|\eta\|_{L^\infty}
  +2|x|\|h\|_\infty\|\dot\eta\|_{L^\infty}
  +\|h\|_\infty\|\dot\eta\|_{L^\infty}^2
  \notag\\
  &\quad\le C\varepsilon^2.
    \label{eq:energy-error}
 \end{align}
 Write $x=r\theta$, where $r=|x|$ and $\theta\in S^{n-1}$.
 It follows that
 \begin{align*}
  d_g(x,0)^2-r^2
  &=\int_0^1g(\gamma)(\dot\gamma,\dot\gamma)\dd s
    -\int_0^1g_{\mathrm E}(sx)(x,x)\dd s\\
  &=\int_0^1|\dot\eta|^2\dd s
    +\int_0^1h(\gamma)(\dot\gamma,\dot\gamma)\dd s\\
  &=\int_0^1h(sx)(x,x)\dd s\pm C\varepsilon^2\\
  &=r\int_0^r h(s\theta)(\theta,\theta)\dd s
    \pm C\varepsilon^2,
 \end{align*}
 where the third equality follows from \eqref{dis-estimate-eta} and \eqref{eq:energy-error}.
 Since $d_g(x,0)=r\pm C\varepsilon$ and $1/2\le r\le3/2$, we obtain
 \[
  d_g(x,0)-r
  =\frac{d_g(x,0)^2-r^2}{d_g(x,0)+r}
  =\frac12\int_0^r h(s\theta)(\theta,\theta)\dd s
    \pm C\varepsilon^2.
 \]
Since $r=1+O(\varepsilon)$ in the boundary annulus,
\[
 \left|
 \int_0^r h(s\theta)(\theta,\theta)\dd s
 -\int_0^1 h(s\theta)(\theta,\theta)\dd s
 \right|
 \le \|h\|_\infty|r-1|
 =O(\varepsilon^2).
\] Consequently, $B_g(0,1)$ lies between the radial domains
 with boundary functions
 \[
  r
  =1-\frac12\int_0^1h(s\theta)(\theta,\theta)\dd s
    \pm C\varepsilon^2.
 \]

 In particular,
 \[
  |B_g(0,1)|
  =\omega_n-\frac12\int_{S^{n-1}}\int_0^1
       h(s\theta)(\theta,\theta)\dd s\dd\theta
       \pm C\varepsilon^2,
  \qquad
  |B_g(0,1)\mathbin{\triangle}B_{g_{\mathrm E}}(0,1)|
  \le C\varepsilon.
 \]
 Finally,
 $\sqrt{\det g}=1+\frac12\tr h\pm C\varepsilon^2$ uniformly.
 Replacing $B_g(0,1)$ by $B_{g_{\mathrm E}}(0,1)$ in the
 integral of $\tr h$ introduces an error of order
 $O(\varepsilon^2)$. This proves \eqref{eq:volume-variation}.
\end{proof}

Next, we express $\mathcal L_n(h)$, up to a quadratic error,
as the pairing of the scalar curvature distribution $R_g$
with a nonpositive radial weight $w_n$.
\begin{lemma}\label{lem:weight}
 Under the assumptions of Lemma~\ref{lem:variation}, define
 \begin{equation}\label{eq:weight}
  w_n(r)=-\frac1{2(n-1)}\int_r^1(s^{2-n}-s)\dd s,
  \qquad 0<r\le1.
 \end{equation}
 We also write $w_n(x)=w_n(|x|)$ and extend it by zero outside
 $B_{g_{\mathrm E}}(0,1)$, assigning an arbitrary value at the origin.
 Then $w_n\in W^{1,1}_0(B_{g_{\mathrm E}}(0,2))$ and
 \begin{equation}\label{eq:weighted-expansion}
  \Vol_g(B_g(0,1))
  =\omega_n+\langle R_g,w_n\rangle+O(\varepsilon^2).
 \end{equation}
 In particular, if $R_g\ge K$ in distributions for a constant $K\in\R$, then
 \begin{equation}\label{eq:unit-upper}
  \Vol_g(B_g(0,1))\le\omega_n-\frac{\omega_n K}{6(n+2)}+C(\varepsilon^2+|\kappa|\varepsilon).
 \end{equation}
\end{lemma}

\begin{proof}
Direct computation gives $w_n\le0$, $w_n(1)=w_n'(1)=0$, and
\begin{equation}\label{eq:weight-hessian}
 D^2w_n-(\Delta w_n)g_{\mathrm E}
 =\frac12g_{\mathrm E}
  -\frac1{2r^{n-1}}\,dr\otimes dr,
 \qquad r=|x|,\quad 0<r<1.
\end{equation}
Moreover, as $r\downarrow0$,
\begin{equation}\label{omega-n-bound}
 w_n(r)=
 \begin{cases}
  O(1),&n=2,\\
  O(1+|\log r|),&n=3,\\
  O(r^{3-n}),&n\ge4,
 \end{cases}
 \qquad
 w_n'(r)=O(r^{2-n}).
\end{equation}
These bounds and $w_n(1)=0$ imply that the zero extension belongs
to $W^{1,1}_0(B_{g_{\mathrm E}}(0,2))$ and can be approximated
in $W^{1,1}$ by nonpositive functions in
$C_c^\infty(B_{g_{\mathrm E}}(0,2))$.

Let
\begin{equation}\label{R-linearization}
 V_{\mathrm{lin}}^k
 :=\partial_a h_{ak}-\partial_k\tr h,
 \qquad
 Lh:=\partial_kV_{\mathrm{lin}}^k
 =\partial_i\partial_jh_{ij}-\Delta\tr h.
\end{equation}
Expanding \eqref{eq:density-divergence} at $g=g_{\mathrm E}$
gives
\[
 |V - V_{\mathrm{lin}}|\le C|h||Dh|.
\]
Moreover, every term in \eqref{eq:density-F} is quadratic
in the first derivatives of $g$, so
\[
 |F|\le C|Dh|^2.
\]
Thus, using \eqref{eq:local-pairing}, we obtain,
in distributions,
\begin{equation}\label{eq:scalar-divergence}
 R_g=Lh+\partial_k(V - V_{\mathrm{lin}})^k+F,
 \qquad
 |V - V_{\mathrm{lin}}|\le C|h||Dh|,
 \qquad
 |F|\le C|Dh|^2.
\end{equation}
These estimates are uniform for $\|h\|_{C^0}$ sufficiently small. We claim that
 \begin{equation}\label{eq:linear-weight}
  \langle Lh,w_n\rangle
  :=-\int_{B_{g_{\mathrm E}}(0,1)}
       V_{\mathrm{lin}}\cdot Dw_n\dd x
  =\mathcal L_n(h),
 \end{equation}
 where $\mathcal L_n(h)$ is defined in \eqref{eq:functional}. First suppose that $h$ is smooth. Set $A_\delta :=B_{g_{\mathrm E}}(0,1)
\setminus\overline{B_{g_{\mathrm E}}(0,\delta)}$, and let $\nu$
be its outward Euclidean unit normal. Integrating by parts gives
\begin{equation}
 \begin{aligned}
  \int_{A_\delta}w_nLh\dd x
  &=-\int_{A_\delta}
       \partial_iw_n(\partial_jh_{ij}-\partial_i\tr h)\dd x
    +\int_{\partial A_\delta}w_nV_{\mathrm{lin}}\cdot\nu\dd S\\
  &=\int_{A_\delta}
       \left\langle D^2w_n-(\Delta w_n)g_{\mathrm E},h
       \right\rangle_{g_{\mathrm E}}\dd x\\
  &\quad+\int_{\partial A_\delta}
       \left(w_nV_{\mathrm{lin}}\cdot\nu
       -h(Dw_n,\nu)+(\tr h)\partial_\nu w_n\right)\dd S\\
  &=\frac12\int_{A_\delta}\tr h\dd x
    -\frac12\int_{S^{n-1}}\int_\delta^1
       h(r\theta)(\theta,\theta)\dd r\dd\theta+o(1).
 \end{aligned}
\end{equation}
Indeed, the outer boundary terms vanish since $w_n(1)=w_n'(1)=0$,
while \eqref{omega-n-bound} bounds the inner boundary terms by
\[
 C\delta^{n-1}\left(
 |w_n(\delta)|\|Dh\|_\infty
 +|w_n'(\delta)|\|h\|_\infty\right)=o(1).
\]
Letting $\delta\downarrow0$ proves \eqref{eq:linear-weight}.
 The tensor on the right-hand side of \eqref{eq:weight-hessian}
 is integrable, and polar integration gives
 \[
\int_{B_{g_{\mathrm E}}(0,1)}
 \left\langle D^2w_n-(\Delta w_n)g_{\mathrm E},h
 \right\rangle_{g_{\mathrm E}}\dd x
  =\frac12\int_{B_{g_{\mathrm E}}(0,1)}\tr h\dd x
   -\frac12\int_{S^{n-1}}\int_0^1
      h(r\theta)(\theta,\theta)\dd r\dd\theta
  =\mathcal L_n(h).
 \]
 For general $C^1$ tensors $h$, smooth approximation in
 $C^1(\overline{B_{g_{\mathrm E}}(0,2)})$ proves
 \eqref{eq:linear-weight}, since $Dw_n\in L^1$ and
 $\mathcal L_n$ is continuous in the $C^0$ norm.

 Since all terms in \eqref{eq:scalar-divergence} define continuous
 functionals on $W^{1,1}_0(B_{g_{\mathrm E}}(0,2))$, we may pair
 this identity with $w_n$. Using \eqref{eq:linear-weight}, we obtain
 \[
  \mathcal L_n(h)
  =\langle R_g,w_n\rangle
   +\int_{B_{g_{\mathrm E}}(0,1)}Dw_n\cdot (V - V_{\mathrm{lin}})\dd x
   -\int_{B_{g_{\mathrm E}}(0,1)}w_n F\dd x.
 \]
 The last two terms have absolute value at most
 \[
  C\varepsilon^2
  \left(
   \|Dw_n\|_{L^1(B_{g_{\mathrm E}}(0,1))}
   +\|w_n\|_{L^1(B_{g_{\mathrm E}}(0,1))}
  \right)
  \le C_n\varepsilon^2.
 \]
 Lemma~\ref{lem:variation} now gives
 \eqref{eq:weighted-expansion}.
 Polar integration gives
 \[
 \int_{B_{g_{\mathrm E}}(0,1)}w_n\dd x
 =-\frac{\omega_n}{2(n-1)}\int_0^1(s^2-s^{n+1})\dd s
 =-\frac{\omega_n}{6(n+2)}.
 \]
 Since $\sqrt{\det g}=1+O(\varepsilon)$,
 \[
 \int_{B_{g_{\mathrm E}}(0,1)}w_n\dd\mu_g
 =-\frac{\omega_n}{6(n+2)}+O(\varepsilon).
 \]
 
 Finally, if $R_g\ge K$, the nonpositive smooth approximations to $w_n$
 and continuity of the pairing imply
 \[
 \langle R_g,w_n\rangle\le
 K\int_{B_{g_{\mathrm E}}(0,1)}w_n\dd\mu_g
 \le-\frac{\omega_n K}{6(n+2)}+C_n|K|\varepsilon,
 \]
 proving \eqref{eq:unit-upper}.
\end{proof}

\begin{remark}
 The choice of the weight in \eqref{eq:weight} can be motivated
 directly from \eqref{eq:functional}. Indeed, polar coordinates give
 \begin{equation}\label{remark-1}
  \mathcal L_n(h)
  =\frac12\int_{B_{g_{\mathrm E}}(0,1)}
   \left\langle
    g_{\mathrm E}-\frac1{r^{n-1}}\,dr\otimes dr,h
   \right\rangle_{g_{\mathrm E}}\dd x.
 \end{equation}
 In view of \eqref{eq:scalar-divergence}, we seek a weight $w_n$
 such that
 $\langle Lh,w_n\rangle=\mathcal L_n(h)$.
 For a weight whose boundary terms vanish as in the proof above,
 integration by parts gives
 \[
  \langle Lh,w_n\rangle
  =-\int_{B_{g_{\mathrm E}}(0,1)}
       V_{\mathrm{lin}}\cdot Dw_n\dd x
  =\int_{B_{g_{\mathrm E}}(0,1)}
    \left\langle
     D^2w_n-(\Delta w_n)g_{\mathrm E},h
    \right\rangle_{g_{\mathrm E}}\dd x.
 \]
 Comparing this expression with \eqref{remark-1} leads to
 \eqref{eq:weight-hessian}.

 For a radial function $w$,
 \[
  D^2w-(\Delta w)g_{\mathrm E}
  =\left(w''-\frac{w'}r\right)dr\otimes dr
   -\left(w''+\frac{n-2}{r}w'\right)g_{\mathrm E}.
 \]
 Matching the coefficients of $dr\otimes dr$ and $g_{\mathrm E}$
 therefore gives
 \[
  w''-\frac{w'}r=-\frac1{2r^{n-1}},
  \qquad
  w''+\frac{n-2}{r}w'=-\frac12.
 \]
 These equations are compatible: subtracting the first from the
 second yields
 \[
  w'(r)=\frac{r^{2-n}-r}{2(n-1)},
 \]
 which also satisfies both equations. Integrating with $w(1)=0$
 gives precisely \eqref{eq:weight}; moreover, $w'(1)=0$ and
 $w\le0$.
\end{remark}

With Lemmas~\ref{lem:variation} and \ref{lem:weight} in hand,
we now prove Theorem~\ref{thm:distributional-balls}.

\begin{proof}[Proof of Theorem~\ref{thm:distributional-balls}]
 Fix $p\in M$. A linear coordinate change followed by a smooth
 quadratic change gives coordinates centered at $p$ with
 \[
  g(0)=g_{\mathrm E},\qquad Dg(0)=0.
 \]
 This normalization uses only the value and first derivatives of
 $g$ at $p$, not an exponential coordinate chart. Choose $r_0>0$ sufficiently small and set
$\rho(r)=\rho_p(r):=\sup_{|x|\le r}|Dg(x)|$ for $0<r<r_0$.
Then $\rho$ is nondecreasing, tends to zero as $r\downarrow0$,
and satisfies
 \begin{equation}\label{eq:single-modulus}
  |Dg(x)|\le\rho(|x|),\qquad
  |g(x)-g_{\mathrm E}|\le |x|\rho(|x|),\qquad |x|<r_0.
 \end{equation}
 The second bound follows by integration along radial segments.
 Let $S_t(x) :=tx$ and define the rescaled metric
 \[
  g_t :=t^{-2}S_t^*g=g_{\mathrm E}+h_t\quad\text{on }B_{g_{\mathrm E}}(0,4),
  \qquad (g_t)_{ij}(x)=g_{ij}(tx).
 \]
 For small $t$,
 \begin{equation}\label{eq:single-rescaled}
  \|h_t\|_{C^1(B_{g_{\mathrm E}}(0,4))}\le Ct\rho(4t).
 \end{equation}

 Under this rescaling, $R_{g_t}\ge t^2K$. Indeed, for every
 nonnegative test function $ u \in C_c^\infty(B_{g_{\mathrm E}}(0,4))$,
 \[
 \langle R_{g_t}, u \rangle
 =t^{2-n}\langle R_g, u \circ S_t^{-1}\rangle
 \ge t^2K\int u \dd\mu_{g_t}.
 \]
 The equality follows by change of variables for smooth metrics
 and by $C^1$ approximation in general.
 Shrink the coordinate neighborhood so that $g\ge\frac12g_{\mathrm E}$.
Any curve starting at $p$ and reaching the coordinate sphere
$\partial B_{g_{\mathrm E}}(0,3t)$ has length at least $\frac{3t}{\sqrt2}>t$. Hence curves of length less than $t$ cannot leave this coordinate
ball. It follows that $D_t$ identifies $B_{g_t}(0,1)$, computed
using curves in $B_{g_{\mathrm E}}(0,4)$, with $B_g(p,t)$. Lemma~\ref{lem:weight} and
 \eqref{eq:single-rescaled} give
 \begin{equation}\label{eq:single-volume}
  \Vol_g(B_g(p,t))
  \le\omega_n t^n-\frac{\omega_n K}{6(n+2)}t^{n+2}+C_n t^{n+2}\bigl(\rho(4t)^2+|K|t\rho(4t)\bigr).
 \end{equation}
 Taking
$\varepsilon_p(t)=C_n\bigl(\rho(4t)^2+|K|t\rho(4t)\bigr)$
proves the first assertion and the volume test.
Under $C^1_{\mathrm{loc}}$ convergence, the first derivatives
are uniformly equicontinuous, and the normalizing coordinate
changes have uniform bounds on a common neighborhood.
Thus the normalized metrics admit a common modulus $\rho$.
The same argument gives a common $t_p$ and $\varepsilon_p$,
proving the last assertion.
\end{proof}

\section{A counterexample below \texorpdfstring{$C^1$}{C1} regularity}
\label{sec:below-c1}

We now prove Theorem~\ref{thm:below-c1}. Let $\sigma$ denote the standard round metric on $\Sph^2$.
We construct the metric in polar coordinates as
\begin{equation}\label{eq:rough-polar-metric}
 g=dr^2+r^2e^{2f(r,\theta)}\sigma,
 \qquad
 f(r,\theta)=r^\beta\phi(\theta)
            +r^{2\beta}\psi(\theta)-ar^{3\beta},
\end{equation}
where $\beta=\sqrt{3}-1$, $a>0$, and
$\phi,\psi\in C^\infty(\Sph^2)$ will be chosen below.
The advantage of this form is that its distance balls centered
at the origin are exactly the coordinate balls:
\begin{equation}\label{eq:rough-distance-balls}
 d_g(0,x)=|x|,
 \qquad
 B_g(0,t)=B_{g_{\mathrm E}}(0,t).
\end{equation}
Indeed, every curve from $0$ to $x$ has length at least $|x|$,
whereas the radial segment has length $|x|$.

We choose $\phi$ and $\psi$ to cancel the first two terms in
the scalar curvature expansion while retaining a positive
leading correction to the ball volume.

\begin{lemma}\label{lem:rough-spherical-data}
 There exist $\phi,\psi\in C^\infty(\Sph^2)$, with
 $\phi\not\equiv0$, such that
 \begin{equation}
  \Delta_\sigma\phi=-6\phi,\quad
  \int_{\Sph^2}\phi\dd\mu_\sigma =0,
  \label{eq:rough-phi}
 \end{equation}
 \begin{equation}
  \bigl(-2\Delta_\sigma-4(2\beta+1)^2\bigr)\psi
    =(20+6\beta^2)\phi^2.
  \label{eq:rough-psi}
 \end{equation}
 Moreover,
 \begin{equation}\label{eq:rough-positive-moment}
  2\int_{\Sph^2}(\psi+\phi^2)\dd\mu_\sigma
  =\frac{\beta^2}{(2\beta+1)^2}
     \int_{\Sph^2}\phi^2\dd\mu_\sigma>0.
 \end{equation}
\end{lemma}

\begin{proof}
 The eigenvalues of $-\Delta_\sigma$ on the standard round sphere
 $\Sph^2$ are $\ell(\ell+1)$, $\ell=0,1,2,\ldots$;
 see e.g. \cite{Chavel1984}.
 Choose a nonzero eigenfunction $\phi$ corresponding to the
 eigenvalue $6$; then \eqref{eq:rough-phi} follows.

 Since $\beta=\sqrt3-1$, the number $2(2\beta+1)^2=26-8\sqrt3$
 is not an eigenvalue of $-\Delta_\sigma$. Consequently,
 the elliptic operator $-2\Delta_\sigma-4(2\beta+1)^2$
 has trivial kernel. By the Fredholm alternative and elliptic
 regularity, there exists a unique $\psi\in C^\infty(\Sph^2)$
 solving \eqref{eq:rough-psi}.
 Integrating over $\Sph^2$ gives
 \[
  \int_{\Sph^2}\psi\dd\mu_\sigma
  =-\frac{20+6\beta^2}{4(2\beta+1)^2}
     \int_{\Sph^2}\phi^2\dd\mu_\sigma.
 \]
Using $\beta^2+2\beta-2=0$, we obtain \eqref{eq:rough-positive-moment}.
\end{proof}

Fix $\phi,\psi$ as in Lemma~\ref{lem:rough-spherical-data}.
We next verify the regularity and distributional curvature
of the metric in \eqref{eq:rough-polar-metric}.

\begin{lemma}\label{lem:rough-nnsc}
 For sufficiently large $a>0$ and sufficiently small $r_0>0$,
 the metric in \eqref{eq:rough-polar-metric}, extended by
 $g(0)=g_{\mathrm E}$, has the regularity asserted in
 Theorem~\ref{thm:below-c1} on $U=B_{g_{\mathrm E}}(0,r_0)$.
 Its distributional scalar curvature is represented by a
 continuous nonnegative function.
\end{lemma}

\begin{proof}
 On the punctured neighborhood, write
$g=dr^2+\gamma_r$, where $\gamma_r=r^2e^{2f}\sigma$.
The induced scalar curvature, second fundamental form, and mean
curvature of the slices are given by
\[
 R_{\gamma_r}
 =\frac{2e^{-2f}}{r^2}(1-\Delta_\sigma f),
 \qquad
 A=\frac12\partial_r\gamma_r
   =\left(\frac1r+f_r\right)\gamma_r,
 \qquad
 H=\tr_{\gamma_r}A=2\left(\frac1r+f_r\right).
\]
Consequently,
\[
 |A|_{\gamma_r}^2=2\left(\frac1r+f_r\right)^2,
 \qquad
 \partial_rH=-\frac2{r^2}+2f_{rr}.
\]
Therefore  we have
\begin{equation}\label{eq:rough-scalar-formula}
 \begin{aligned}
  R_g
  &=R_{\gamma_r}-|A|_{\gamma_r}^2-H^2-2\partial_rH\\
  &=\frac{2e^{-2f}}{r^2}(1-\Delta_\sigma f)
    -\frac2{r^2}-\frac{12}{r}f_r-6f_r^2-4f_{rr}.
 \end{aligned}
\end{equation}
To compute the leading coefficients, expand
\eqref{eq:rough-scalar-formula} as
\[
 \begin{aligned}
 r^2R_g
 ={}&-2\Delta_\sigma f-4f-12rf_r-4r^2f_{rr}\\
    &+4f^2+4f\Delta_\sigma f-6r^2f_r^2\\
    &-\frac83f^3-4f^2\Delta_\sigma f
      +O_a(r^{4\beta}).
 \end{aligned}
\]
Substituting \eqref{eq:rough-polar-metric} and using
$\Delta_\sigma\phi=-6\phi$, we obtain
 \begin{align}
  r^2R_g
  ={}&-4(\beta^2+2\beta-2)\phi r^\beta \notag\\
   &+\bigl[-2\Delta_\sigma\psi-4(2\beta+1)^2\psi
             -(20+6\beta^2)\phi^2\bigr]r^{2\beta}\notag\\
   &+\bigl[4(3\beta+1)^2a+Q\bigr]r^{3\beta}
       +O_a(r^{4\beta}),
  \label{eq:rough-scalar-expansion}
 \end{align}
 uniformly in $\theta$, where
 \[
  Q=4\phi\Delta_\sigma\psi
       -(16+24\beta^2)\phi\psi+\frac{64}{3}\phi^3.
 \]
 The first two coefficients vanish by the choice of $\beta$
 and Lemma~\ref{lem:rough-spherical-data}. Choose $a$ so that
 \[
  4(3\beta+1)^2a>\|Q\|_\infty+1.
 \]
 After decreasing $r_0$, we have $R_g>0$ for $0<r<r_0$.
 Since $3\beta-2>0$, this scalar curvature extends continuously
 to the origin by $R(0)=0$.

 In Cartesian coordinates, with $\theta=x/|x|$,
 \[
  g_{ij}=\delta_{ij}
       +(e^{2f}-1)(\delta_{ij}-\theta_i\theta_j).
 \]
 Thus
 \begin{equation}\label{eq:rough-metric-bounds}
  |g-g_{\mathrm E}|\le Cr^\beta,
  \qquad |Dg|\le Cr^{\beta-1}.
 \end{equation}
 These bounds imply $g\in C^{0,\beta}(U)$. Indeed, for points whose
 separation is comparable to their distances from the origin,
 use the first bound; otherwise apply the second along the
 segment joining them. The derivatives on the punctured
 neighborhood are also weak derivatives across the origin,
 since the boundary terms in integration by parts are
 $O(\delta^2)$. Consequently,
 \[
  g\in W^{1,p}_{\mathrm{loc}}(U)
  \quad\text{for}\quad
  p<\frac3{1-\beta}=6+3\sqrt3.
 \]

 Finally, we check the distributional scalar curvature at the
origin. In Cartesian coordinates,
\[
 J=e^{2f},
 \qquad
 Jg^{ij}=\delta_{ij}+(e^{2f}-1)\theta_i\theta_j.
\]
Using
\[
 g^{ab}\Gamma^k_{ab}
 =-\frac1J\partial_i(Jg^{ik}),
 \qquad
 \Gamma^b_{ab}=\partial_a\log J,
\]
we can rewrite \eqref{eq:density-divergence} as
\[
 V^k=-\partial_i(Jg^{ik})-g^{ik}\partial_iJ.
\]
This gives
\begin{equation}\label{eq:rough-V-F}
 V=-\left(4e^{2f}f_r+\frac{2(e^{2f}-1)}r\right)\theta
      -\frac2r\nabla_\sigma f,
 \qquad
 F=2e^{2f}f_r^2,
\end{equation}
where $\nabla_\sigma f$ is viewed as a Euclidean vector tangent
to the unit sphere. Indeed, differentiating the expression for $V$ gives
\begin{equation}
 \operatorname{div}V
 =-4e^{2f}f_{rr}-8e^{2f}f_r^2
     -\frac{12e^{2f}}r f_r-\frac{2(e^{2f}-1)}{r^2}
     -\frac2{r^2}\Delta_\sigma f.
\end{equation}
Comparison with \eqref{eq:rough-scalar-formula} yields
$F=JR_g-\operatorname{div}V=2e^{2f}f_r^2$. Thus
\[
 |V|=O(r^{\beta-1}),\qquad |F|=O(r^{2\beta-2}),
\]
so both are locally integrable. For $u \in C_c^\infty(U)$, set
$U_\delta=U\setminus\overline{B_{g_{\mathrm E}}(0,\delta)}$,
and let $\nu_\delta$ denote the outward unit normal of
$B_{g_{\mathrm E}}(0,\delta)$. Since
$JR=\operatorname{div}V+F$ on $U_\delta$, integration by parts gives
\begin{equation}
 \begin{aligned}
  \int_{U_\delta}R u \dd\mu_g
  &=\int_{U_\delta} u (\operatorname{div}V+F)\dd x\\
  &=-\int_{U_\delta}V\cdot D u \dd x
    +\int_{U_\delta}F u \dd x
    -\int_{\partial B_{g_{\mathrm E}}(0,\delta)}
        u V\cdot\nu_\delta\dd S.
 \end{aligned}
\end{equation}
There is no outer boundary contribution because $u$ has
compact support in $U$. The inner boundary term tends to zero, and hence
\[
 \left|
  \int_{\partial B_{g_{\mathrm E}}(0,\delta)}
       u V\cdot\nu_\delta\dd S
 \right|
 \le \|u \|_\infty
      \sup_{|x|=\delta}|V(x)|\,4\pi\delta^2
 \le C\|u \|_\infty\delta^{\beta+1}
 \longrightarrow0.
\]
Since $V,F\in L^1_{\mathrm{loc}}(U)$ and $R$ is continuous,
letting $\delta\downarrow0$ yields
\[
 \langle R_g, u\rangle
 =-\int_U V\cdot D u \dd x+\int_U F u \dd x
 =\int_U R u\dd\mu_g.
\]
Thus there is no additional curvature distribution supported
at the origin, and $R_g\ge0$.
\end{proof}

\begin{proof}[Proof of Theorem~\ref{thm:below-c1}]
 Choose the metric given by Lemma~\ref{lem:rough-nnsc}.
 It remains to compute its ball volumes.
 By \eqref{eq:rough-distance-balls},
 \[
  \Vol_g(B_g(0,t))
  =\int_0^t r^2\int_{\Sph^2}e^{2f(r,\theta)}
       \dd\mu_\sigma\dd r.
 \]
 Uniformly in $\theta$,
 \[
  e^{2f}
  =1+2r^\beta\phi
     +2r^{2\beta}(\psi+\phi^2)+O(r^{3\beta}).
 \]
 Applying \eqref{eq:rough-phi} and
 \eqref{eq:rough-positive-moment}, we obtain
 \[
  \Vol_g(B_g(0,t))
  =\frac{4\pi}{3}t^3
   +\frac{\beta^2 t^{3+2\beta}}
          {(2\beta+1)^2(3+2\beta)}
      \int_{\Sph^2}\phi^2\dd\mu_\sigma
   +O(t^{3+3\beta}).
 \]
 The coefficient of $t^{3+2\beta}$ is strictly positive.
 Since $2\beta-2<0$ and $3\beta-2>0$, division by $t^5$
 proves the asserted divergence.
\end{proof}

\begin{remark}
 The preceding example also shows that the $C^1_{\mathrm{loc}}$
 convergence in Theorem~\ref{thm:distributional-balls} cannot
 in general be replaced by $C^0_{\mathrm{loc}}$ convergence.
 Indeed, fix $3<p<6+3\sqrt3$. After shrinking the neighborhood
 of the origin, extend $g$ to a continuous $W^{1,p}$ metric
 $\widehat g$ on a closed three-manifold, smooth away from the
 origin and agreeing with $g$ near it. Since the scalar curvature
 is nonnegative near the origin and smooth elsewhere, there is
 a constant $K_0\le0$ such that $R_{\widehat g}\ge K_0$ in the
 distributional sense.

 By the scalar curvature preserving regularization of
 Jiang--Sheng--Zhang \cite[Theorem~1.1]{JiangShengZhang2023},
 using the Ricci--DeTurck gauge, there exist smooth metrics
 $\widehat g_j$ such that
 \[
  \widehat g_j\longrightarrow\widehat g
  \quad\text{in }C^0,
  \qquad R_{\widehat g_j}\ge K_0.
 \]
 Nevertheless, the local ball volume calculation above gives
 \[
  \lim_{t\downarrow0}
  \frac{\Vol_{\widehat g}(B_{\widehat g}(0,t))
        -\omega_3t^3}{t^5}=+\infty.
 \]
 Thus the limit fails the volume limit test for $K_0$, even though
 the smooth approximating metrics have scalar curvature
 bounded below by $K_0$.
\end{remark}

\section{An example for fixed tangent regions}\label{sec:example}

We present the reflected-sector construction from
\cite[Section~6]{WangWang2026} in detail.
The volume calculation uses the Ricci tensors of two smooth
local models.

\subsection{Curvature formulas}

Choose $p\in\Sph^2$ and coordinates $(r,\theta)$ centered at $p$
such that the standard round metric is
\[
 h=dr^2+\cos^2r\,d\theta^2.
\]

\begin{lemma}\label{lem:scalar-formula}
 Let $h$ be the unit round metric on $\Sph^2$ and let
 $\varphi$ be a smooth positive function. The metric
 $g=\varphi^{-2}h+\varphi^2d\xi^2$ satisfies
 \begin{equation}\label{eq:scalar-formula}
  R_g=2\varphi^2-2|\nabla_h\varphi|_h^2.
 \end{equation}
In coordinates $h=dr^2+\cos^2r\,d\theta^2$,
 \begin{equation}\label{eq:mixed-ricci}
  \Ric_g(\partial_r,\partial_\theta)
  =-\frac{\varphi_{r\theta}}{\varphi}
   -2\frac{\varphi_r\varphi_\theta}{\varphi^2}
   -\tan r\,\frac{\varphi_\theta}{\varphi}.
 \end{equation}
 Moreover,
 $\Ric_g(\partial_r,\partial_\xi)
 =\Ric_g(\partial_\theta,\partial_\xi)=0$.
\end{lemma}
\begin{proof}
 Put $\bar h=\varphi^{-2}h$. For $a,b\in\{r,\theta\}$,
 the warped product formulas give
 \begin{equation}\label{Ric-computation}
  \begin{aligned}
  \Ric_g(\partial_a,\partial_b)
  &=\Ric_{\bar h}(\partial_a,\partial_b)
   -\varphi^{-1}\Hess_{\bar h}\varphi(\partial_a,\partial_b)\\
  \Ric_g(\partial_\xi,\partial_\xi)
  &=-\varphi\Delta_{\bar h}\varphi
  =-\varphi^3\Delta_h\varphi\\
\Ric_g(\partial_a,\partial_\xi)&=0
  \end{aligned}
 \end{equation}
 The conformal transformation formulas give
 \begin{equation}
 \begin{aligned}\label{Ric-conformal-computation}
    \Ric_{\bar h}(\partial_a,\partial_b)
  &=\left(
    1+\frac{\Delta_h\varphi}{\varphi}
    -\frac{|\nabla_h\varphi|_h^2}{\varphi^2}
   \right)h_{ab}\\
  \Hess_{\bar h}\varphi(\partial_a,\partial_b)
  &=\Hess_h\varphi(\partial_a,\partial_b)
   +\frac{2\varphi_a\varphi_b}{\varphi}
   -\frac{|\nabla_h\varphi|_h^2}{\varphi}h_{ab}.
 \end{aligned}
 \end{equation}
 Substituting \eqref{Ric-conformal-computation} into \eqref{Ric-computation} and cancelling the gradient terms proportional to
 $h_{ab}$, we obtain
 \[
  \Ric_g(\partial_a,\partial_b)
  =\left(1+\frac{\Delta_h\varphi}{\varphi}\right)h_{ab}
   -\frac{\Hess_h\varphi(\partial_a,\partial_b)}{\varphi}
   -2\frac{\varphi_a\varphi_b}{\varphi^2}.
 \]
 Taking the trace with respect to $g$ yields
 \[
  R_g
  =2\varphi^2+\varphi\Delta_h\varphi
   -2|\nabla_h\varphi|_h^2-\varphi\Delta_h\varphi
  =2\varphi^2-2|\nabla_h\varphi|_h^2.
 \]
 Since $h_{r\theta}=0$ and
 \[
  \Hess_h\varphi(\partial_r,\partial_\theta)
  =\varphi_{r\theta}+\tan r\,\varphi_\theta,
 \]
 we also obtain
 \[
  \Ric_g(\partial_r,\partial_\theta)
  =-\frac{\varphi_{r\theta}}{\varphi}
   -2\frac{\varphi_r\varphi_\theta}{\varphi^2}
   -\tan r\,\frac{\varphi_\theta}{\varphi}.
 \]
\end{proof}

\subsection{Construction of the metric and region}

Fix the target exponent $0<\alpha<1$, and choose
$\alpha<\gamma<1$. Let $b\in C^\infty(\R)$ satisfy
\[
 0\le b\le1,\qquad b(s)=0\ \text{for }s\le\tfrac12,
 \qquad b(s)=1\ \text{for }s\ge1.
\]
On a small square centered at the origin, define
\begin{equation}\label{eq:profile}
 \psi(r,\theta)=
 \begin{cases}
  |r\theta|\,
  b\!\left(\dfrac{|\theta|}{|r|^{1/\gamma}}\right)
  b\!\left(\dfrac{|r|}{|\theta|^{1/\gamma}}\right),
       &r\theta\ne0,\\[2mm]
  0,&r\theta=0.
 \end{cases}
\end{equation}
Thus $\psi=|r\theta|$ where
$|\theta|\ge|r|^{1/\gamma}$ and $|r|\ge|\theta|^{1/\gamma}$,
and $\psi$ vanishes in a neighborhood of each nonzero point of
either coordinate axis.

\begin{lemma}\label{lem:profile}
 The function $\psi$ is smooth away from the origin and belongs to
 $C^{1,\gamma}$ near the origin. Moreover,
 \begin{equation}\label{eq:profile-first}
  |\psi|\le|r\theta|,\qquad
  |\psi_r|\le C|\theta|,\qquad
  |\psi_\theta|\le C|r|,
 \end{equation}
 and, away from the coordinate axes,
 \begin{equation}\label{eq:profile-second}
  |\psi_{rr}|\le C|r|^{\gamma-1},\qquad
  |\psi_{r\theta}|\le C,\qquad
  |\psi_{\theta\theta}|\le C|\theta|^{\gamma-1}.
 \end{equation}
\end{lemma}

\begin{proof}
In a fixed quadrant put
 \[
  A=\frac{|\theta|}{|r|^{1/\gamma}},\qquad
  B=\frac{|r|}{|\theta|^{1/\gamma}},\qquad
  \sigma=\operatorname{sgn}(r\theta).
 \]
 Then
 \[
  \psi_r=\sigma\theta
    \left[b(A)b(B)-\gamma^{-1}Ab'(A)b(B)+Bb(A)b'(B)\right].
 \]
 Since $sb'(s)$ is bounded, this and the corresponding formula
 with $r,\theta$ interchanged prove \eqref{eq:profile-first}.
 Differentiating once more gives
 \[
  |\psi_{rr}|\le C\frac{|\theta|}{|r|},\qquad
  |\psi_{r\theta}|\le C,\qquad
  |\psi_{\theta\theta}|\le C\frac{|r|}{|\theta|}.
 \]
 Wherever the relevant derivatives can be nonzero, $A,B>1/2$,
 and therefore $|\theta|\le(2|r|)^\gamma$ and
 $|r|\le(2|\theta|)^\gamma$. This proves
 \eqref{eq:profile-second}.

 The first estimates show that $\psi$ is differentiable at the
 origin, with $D\psi(0)=0$, and that its gradient is continuous
 there. To prove H\"older continuity, the triangle inequality and
\eqref{eq:profile-second} give
\begin{align*}
 |\psi_r(r_1,\theta_1)-\psi_r(r_2,\theta_2)|
 &\le
 |\psi_r(r_1,\theta_1)-\psi_r(r_2,\theta_1)|
 +|\psi_r(r_2,\theta_1)-\psi_r(r_2,\theta_2)|\\
 &\le
 C\int_{\min\{r_1,r_2\}}^{\max\{r_1,r_2\}}
 |s|^{\gamma-1}\dd s
 +C|\theta_1-\theta_2|\\
 &\le C\bigl(|r_1-r_2|^\gamma
             +|\theta_1-\theta_2|^\gamma\bigr),
\end{align*}
where we used $
 \int_a^b |s|^{\gamma-1}\dd s
 \le \frac{2}{\gamma}|b-a|^\gamma$ for $a<b$
and bounded the linear term by its $\gamma$-power on a fixed
small square. If a segment crosses an axis, we split the integral
there and use continuity. The same argument applies to
$\psi_\theta$. Hence $\psi\in C^{1,\gamma}$.
\end{proof}

Let $\chi\in C_c^\infty((-2,2)^2)$ satisfy $0\le\chi\le1$ and
$\chi=1$ on $[-1,1]^2$, and put
\[
 \chi_\rho(r,\theta)=\chi(r/\rho,\theta/\rho).
\]
For $k>0$, to be chosen below, define
\begin{equation}\label{eq:global-phi}
 \varphi=1+k\chi_\rho\psi,
 \qquad g=\varphi^{-2}h+\varphi^2d\xi^2,
\end{equation}
extending $\chi_\rho\psi$ by zero outside the coordinate chart.
Lemma~\ref{lem:profile} gives
\[
 \|\chi_\rho\psi\|_\infty\le4\rho^2,
 \qquad
 \|\nabla_h(\chi_\rho\psi)\|_\infty\le C\rho.
\]
Choose $\rho>0$ sufficiently small that
\begin{equation}\label{eq:phi-bounds}
 |\varphi-1|<\frac14,
 \qquad |\nabla_h\varphi|_h<\frac12.
\end{equation}
Then $g\in C^{1,\gamma}$, is smooth off
$\Gamma=\{p\}\times\Sph^1$, and agrees with the product metric
outside the support of $\chi_\rho$.

Fix
\begin{equation}\label{eq:omega}
 \Omega=\{0\}\cup
 \left\{v=(x,y,z)\in\R^3:
 0<|v|<1,\quad
 \min\{|x|,|y|,|z|\}>\frac14|v|\right\}.
\end{equation}
This set consists of eight truncated conical sectors around the
directions $(\pm1,\pm1,\pm1)$, together with the origin.
It is bounded, measurable, star-shaped, and has positive volume.
Since independent sign changes and coordinate permutations preserve
$\Omega$,
\begin{equation}\label{eq:omega-moments}
 \int_\Omega x^2\dd v
 =\int_\Omega y^2\dd v
 =\int_\Omega z^2\dd v=:m>0,
 \qquad
 \int_\Omega xy\dd v
 =\int_\Omega xz\dd v
 =\int_\Omega yz\dd v=0.
\end{equation}
Set
\begin{equation}\label{eq:A}
 A=\int_\Omega|xy|\dd v>0.
\end{equation}
Both $m$ and $A$ are independent of $\alpha,k,\rho$.

Let $\Gamma=\{p\}\times\Sph^1$ as in Theorem~\ref{thm:example}. At each $q\in\Gamma$,
identify $T_qM$ with $\R^3$ using the orthonormal coordinate
frame $(\partial_r,\partial_\theta,\partial_\xi)$.
Near $q$, the cutoff equals one. For $\sigma\in\{+1,-1\}$,
consider the smooth local metrics
\begin{equation}\label{eq:model}
 g^\sigma
 =(1+\sigma kr\theta)^{-2}(dr^2+\cos^2r\,d\theta^2)
 +(1+\sigma kr\theta)^2d\xi^2.
\end{equation}
They agree with $g$ in the respective regions
\[
 P^\sigma=
 \left\{\operatorname{sgn}(r\theta)=\sigma,\quad
 |\theta|>|r|^{1/\gamma},\quad
 |r|>|\theta|^{1/\gamma}\right\}
\]
inside this neighborhood. Write
\[
 \Omega^\sigma
 =\{v\in\Omega:\operatorname{sgn}(xy)=\sigma\}.
\]

\begin{lemma}\label{lem:admissible}
For every $q\in\Gamma$ and all sufficiently small $t>0$,
the exponential map $\exp_q^g$ is well-defined and injective
on $t\Omega$, its image is measurable, and
\begin{equation}\label{eq:model-agreement}
 \exp_q^g(t\Omega^\sigma)
 =\exp_q^{g^\sigma}(t\Omega^\sigma)\subset P^\sigma.
\end{equation}
\end{lemma}

\begin{proof}
 Consider unit directions $v=(a,b,c)$ in the closure of the
 directional set in \eqref{eq:omega}, so that $|a|,|b|\ge1/4$.
 The Christoffel symbols of $g^\sigma$ vanish at $q$. Hence its
 geodesics with initial velocity $v$ satisfy
 \[
  r(s)=as+O(s^3),\qquad
  \theta(s)=bs+O(s^3),\qquad
  \xi(s)=cs+O(s^3),
 \]
 uniformly over these directions. Taking $\sigma=\operatorname{sgn}(ab)$ and using
$1/4\le |a|,|b|\le1$, we obtain, for all sufficiently small
$s>0$, uniformly over these directions,
\[
 \operatorname{sgn}(r(s)\theta(s))=\sigma,
 \qquad
 \frac{|\theta(s)|}{|r(s)|^{1/\gamma}},\quad
 \frac{|r(s)|}{|\theta(s)|^{1/\gamma}}
 \ge C_0 s^{1-1/\gamma}>1.
\]
 Thus these geodesics stay in $P^\sigma$ and also solve the
 geodesic equation for $g$.

We next prove uniqueness of the $g$-geodesic starting at $q$
with initial velocity $v$. Any such geodesic satisfies
 \[
  r(s)=as+o(s),\qquad \theta(s)=bs+o(s),
 \]
 and therefore enters $P^\sigma$ immediately. It consequently
 solves the smooth comparison equation with the same initial
 data at $q$. Uniqueness for the smooth comparison equation proves
\eqref{eq:model-agreement}. The case of zero initial velocity
is trivial.

 For small $t$, each smooth exponential map is a diffeomorphism
on a neighborhood containing $t\Omega^\sigma$, so its image
is measurable. The two images lie in the disjoint regions
$P^+$ and $P^-$. Together with $\exp_q^g(0)=q$, this shows
that $\exp_q^g$ is injective on $t\Omega$ and has a measurable
image.
\end{proof}

We now prove Theorem~\ref{thm:example}.

\begin{proof}
 Given $\alpha$ and $K$, choose $\gamma\in(\alpha,1)$ and
 $k>(m+3K)/A$. Choose $\rho$ so that \eqref{eq:phi-bounds}
 holds with strict inequalities, and construct $\varphi$ and $g$
 as above.

 Applying Lemma~\ref{lem:scalar-formula} to the smooth model
$g^\sigma$ at $q$ for each sector, we obtain
 \begin{equation}\label{eq:ricci-model}
  \Ric_{g^\sigma}(q)=
  \begin{pmatrix}
   1&-\sigma k&0\\
   -\sigma k&1&0\\
   0&0&0
  \end{pmatrix}.
 \end{equation}
 By Lemma~\ref{lem:admissible}, the volume on each tangent sector
 can be computed using its smooth comparison metric. Applying the first line of \eqref{omega-volume-expansion}
to each tangent sector gives
 \[
  \Vol_g(\exp_q^g(t\Omega^\sigma))
  =t^3|\Omega^\sigma|
  -\frac{t^5}{6}\int_{\Omega^\sigma}
  (x^2+y^2-2\sigma kxy)\dd v+O(t^6).
 \]
 Since the images are disjoint and $\sigma xy=|xy|$ on
 $\Omega^\sigma$, summing gives
 \begin{equation}\label{eq:final-sector-volume}
  \Vol_g(\exp_q^g(t\Omega))-t^3|\Omega|
  =\frac{kA-m}{3}t^5+O(t^6).
 \end{equation}
 Our choice of $k$ proves \eqref{eq:bad-volume} at every
 $q\in\Gamma$, since translation in $\xi$ is an isometry.

 Since $\varphi\in C^{1,\gamma}$ and $\gamma>\alpha$, we can
 choose smooth positive functions $\varphi_j\to\varphi$ in
 $C^{1,\alpha}$ satisfying
 \[
  |\varphi_j-1|\le\frac14,
  \qquad |\nabla_h\varphi_j|_h\le\frac12.
 \]
 By Lemma~\ref{lem:scalar-formula}, the smooth metrics
 $g_j=\varphi_j^{-2}h+\varphi_j^2d\xi^2$ satisfy $R_{g_j}\ge5/8$
 and converge to $g$ in $C^{1,\alpha}$. Thus the $\Omega$-volume
 test fails under this convergence, even with a uniform positive
 scalar curvature lower bound. Moreover, this convergence implies
 $R_g\ge5/8$ in the distributional sense.
\end{proof}

By Theorem~\ref{thm:distributional-balls}, this metric still
satisfies the distance ball volume limit test. The failure above
concerns exponential images of fixed tangent regions.


\bibliographystyle{plain}
\bibliography{volume-tests-ref}
\end{document}